\documentclass[letterpaper, 10pt, conference]{ieeeconf}      %

\IEEEoverridecommandlockouts                              %
\usepackage[T1]{fontenc}
\usepackage{indentfirst}
\usepackage{amssymb}
\usepackage{graphicx}
\usepackage{subfig}
\usepackage{mathptmx}
\usepackage{amsmath}

\newtheorem{de}{Definition}[section]
 
\newtheorem{lemma}{Lemma}[section]
\newtheorem{theo}{Theorem}[section]

\newcommand{\dist}{\operatorname{dist}}

\newcommand{\Tr}{\operatorname{Tr}}

\newcommand{\R}{\mathbb{R}}
\newcommand{\rmax}{\R_{\max}}
\newcommand{\deltat}{\tau}

\newcommand{\Span}{\operatorname{span}}
\newcommand{\Spanb}{\overline{\operatorname{span}}}
\newcommand{\cB}{{\cal{B}}}
\newcommand{\cS}{{\cal{S}}}
\newcommand{\cW}{{\cal{W}}}
\newcommand{\cZ}{{\cal{Z}}}
\newcommand{\setomega}{{W}}

\title{\LARGE \bf
Curse of dimensionality reduction in max-plus based approximation methods: theoretical estimates and improved pruning algorithms
}

\author{ \parbox{2 in}{\centering St\'ephane Gaubert$^*$
         \thanks{$^*$These authors were partially supported by the Arpege program of the French National Agency of Research (ANR), project ``ASOPT'', number ANR-08-SEGI-005 , by the Digiteo project DIM08 ``PASO'' number 3389.}
\\
          INRIA and CMAP\\
          \'Ecole Polytechnique\\ 
         91128 Palaiseau C\'edex, France\\
         {\tt\small Stephane.Gaubert@inria.fr}\\ \
}
         \hspace*{ 0.0 in}
         \parbox{2.5 in}{ \centering William McEneaney$^{**}$
         \thanks{$^{**}$ This author acknowledges support from AFOSR.}
\thanks{Prepared for the 50th IEEE Conference on Decision and Control and European Control Conference (CDC-ECC 2011), \copyright 2011 IEEE}
\\
 Department of Mechanical and Aerospace Engineering\\
 University of California, San Diego\\
 La Jolla, CA 92093-0411, USA
         {\tt\small wmceneaney@eng.ucsd.edu}
}
         \hspace*{ 0.0 in}
         \parbox{2 in}{ \centering Zheng Qu$^*$%
\\
         CMAP, INRIA and Fudan University\\
         \'Ecole Polytechnique\\
         91128 Palaiseau C\'edex, France\\
         {\tt\small zheng.qu@polytechnique.edu}
\\ \ 
}
}

\newcommand{\bx}{\mathbf{x}}
\newcommand{\bu}{\mathbf{u}}
\begin{document}

\maketitle
\thispagestyle{empty}
\pagestyle{empty}

\begin{abstract}
Max-plus based methods have been recently developed to approximate the
value function of possibly high dimensional optimal control
problems. A critical step of these methods consists in approximating a
function by a supremum of a small number of functions (max-plus
``basis functions'') taken from a prescribed dictionary.  We study
several variants of this approximation problem, which we show to be
continuous versions of the facility location and $k$-center
combinatorial optimization problems, in which the connection costs
arise from a Bregman distance. We give theoretical error estimates,
quantifying the number of basis functions needed to reach a prescribed
accuracy. We derive from our approach a refinement of the curse of
dimensionality free method introduced previously by McEneaney, with a
higher accuracy for a comparable computational cost.
\end{abstract}

\section{Introduction}
Dynamic programming is one of the main approaches to optimal control.
It leads to solving Hamilton-Jacobi-Bellman (HJB) partial differential equations. 
Several techniques have been proposed in the literature to solve
this problem. We mention, for example, finite difference schemes
and the method of the vanishing viscosity \cite{crandall-lions},
the antidiffusive schemes for advection \cite{zidani-bokanowski},
the so-called discrete dynamic programming method or
semi-Lagrangian method \cite{capuzzodolcetta,%
falcone,%
falcone-ferretti,%
carlini-falcone-ferretti}.
Unlike alternative approaches based on the maximum principles
or on direct methods, dynamic programming based methods are guaranteed to give
the global optimum of the problem. However, they suffer from
the curse of dimensionality, meaning that the execution time
grows exponentially with the dimension of the state space. 

Recently, a new class of methods has been developed after
the work of Fleming and McEneaney~\cite{a5}, see in particular~\cite{curseofdim,a6,a7}. These methods all rely on max-plus algebra. 
Their common idea is to approximate the value function by a supremum
of finitely many ``basis functions''.  They exploit the max-plus
linearity of the Lax-Oleinik semi-group (evolution semi-group of the HJB partial differential equations associated to a deterministic optimal control problem).
One of these methods, developed by Akian and al.~\cite{a6}, was shown
to be a max-plus analogue of the (Petrov-Galerkin)
finite element method. In particular,
the global error of the method can be estimated in terms of certain
elementary projection errors, as in the case of usual finite elements.

Among the max-plus methods, the curse of dimensionality reduction
of McEneaney~\cite{curseofdim} (see also~\cite{mccomplex,a7,eneaneyphys}) appears to be of special interest.
In its original form, it applies to an optimal switching problem involving
$m$ linear quadratic models: it approximates
the solution by a supremum of quadratic functions
which are obtained by solving Riccati equations. 
The theoretical analysis of the method~\cite{mccomplex} shows that
the growth of the execution time is only polynomial as the dimension grows,
keeping all the other parameters fixed. However, the bound
of~\cite{mccomplex} still grows
exponentially as the required accuracy tends to zero, hence
the curse of dimensionality is replaced by a ``curse of complexity''~\cite{a7}.
However, the complexity of the method can be considerably reduced
in practice by incorporating a pruning algorithm, which eliminates
on the fly the redundant basis functions produced by the algorithm. 
In this way, high dimensional instances (with state dimensions 
from 6 to 15) inaccessible by other methods 
could be solved~\cite{a7,eneaneyphys}. 

This raises the question to understand why, and to what extent,
max-plus techniques can attenuate the curse of dimensionality:
this is the object of the present paper. 

After a brief review of max-plus based methods (Section~\ref{section-review})
we establish in Section~\ref{sec-curseunavoidable}
a negative result, concerning the family of max-plus methods
based on $c$-semiconvex transforms, developed by Fleming and McEneaney~\cite{a5}
and Akian et al.~\cite{a6}. In these methods, the function is approximated
by a supremum of quadratic forms all of which have the same hessian.
Then, Theorem~\ref{theo-unavoidable} below shows that the number of max-plus basis functions
necessary to reach an accuracy of $\epsilon$ is at least of order
$\epsilon^{-d/2}$, meaning
that the curse of dimensionality is inherent to all of these methods (the
order $\epsilon^{-d/2}$ is optimal,
it is reached in particular by the max-plus finite elements of~\cite{a6},
with $P_2$ finite elements and $P_1$ or $P_2$ test functions~\cite{theseasma}). The proof, which is sketched in Section~\ref{sec-proof}, relies on results concerning the approximation of smooth convex bodies~\cite{a2,b2}.

However, this theoretical negative result is contrasted by the experimental
efficiency of pruning in the dimensionality free method of~\cite{curseofdim},
which often gives approximations of an acceptable accuracy for
a modest amount of basis functions.
Therefore, we focus our attention on the algorithmic aspects of the pruning problem in the rest of the paper.
In Section~\ref{sec3}, we present a {\em primal} variant of the
method,  which avoids the use of dual representations:
in the absence of pruning,
it is equivalent to the original method, but we 
shall see that it leads to a more efficient pruning.
Next, we show in Section~\ref{sec4}
that the optimal pruning problem can be formulated as a continuous
version of the $k$-median or $k$-center problem, depending on the choice
of the norm. The discrete versions of these problems are NP-hard.
Hence, we propose several heuristics (combining facility location
heuristics and Shor SDP relaxation scheme).
Experimental results are given in Section~\ref{sec-exp}. They show
that by combining the primal version of the method with improved
pruning algorithms, a higher accuracy is reached for a similar
running time, by comparison with~\cite{curseofdim,a7}.

\section{Max-plus numerical methods to solve optimal control problems}\label{section-review}
\subsection{The Lax-Oleinik semi-group}

We consider the optimal control problem
\begin{equation}\label{dyuj1}
 v(x,T):=\sup \int_0^T \ell(\bx(s),\bu(s)) ds +\phi (\bx(T)) \enspace;
\end{equation}
\begin{equation}
\dot \bx(s)=f(\bx(s),\bu(s)),\quad \bx(0)=x,\quad \bx(s)\in X, \bu(s)\in U
\enspace .\label{e-dynsys}
\end{equation}
Here, $X\subset \mathbb{R}^d$ is the set of states, $U\subset \mathbb{R}^m$ is the set of actions, $T$ denotes the horizon, the initial condition $x\in X$,
the {Lagrangian} $\ell:X\times U \to \R$, the terminal reward $\phi:\R\to\R$, and the dynamics $f:X\times U\to\R^d$ are given.
The supremum is taken over all the control functions $\bu$ and system trajectories $\bx$ satisfying~\eqref{e-dynsys}, and $v$ is the \emph{value function}. We will assume here for simplicity that the set $X$ is invariant by the dynamics~\eqref{e-dynsys} for all choices of the control function $\bu$.
Under certain regularity assumptions, it is known that $v(x,t)$ is the solution of the Hamilton-Jacobi equation:
\begin{equation}\label{guih}
 -\frac{\partial v}{\partial t}+H(x,\frac{\partial v}{\partial x})=0,\qquad \forall(x,t)\in X\times (0,T],
\end{equation}
with initial condition:
\begin{equation}
 v(x,0)=\phi(x),\quad \forall x\in X \enspace .
\end{equation}
Let $(S_T)_{T\geq 0}$ be the \emph{Lax-Oleinik} semi-group, i.e., the \textit{evolution semi-group} of the Hamilton-Jacobi equation.
For every horizon $T$, $S_T$ is a map which associates to the terminal
reward $\phi$ the value function $S_T[\phi]:=v(x,T)$ on horizon $T$.
By \emph{semi-group}, we mean that $S_{t+s}=S_t\circ S_s$ for all $t,s\geq 0$.
Recall that the \textit{max-plus semiring}, $\R_{max}$, is the set $\R\cup \{-\infty\}$, equipped with the addition $(a,b)\mapsto \max(a,b)$ and the multiplication $(a, b)\mapsto a+b$. 
For all maps $f,g$ from $X$ to $\R_{max}$ and $\lambda \in \R_{max}$, we denote
by $f\vee g$ the map such that $(f\vee g) (x)=\max(f(x), g(x))$ and 
by $\lambda+f$ the map such that $(\lambda +f) (x)=\lambda + f(x)$.
It is known that the semi-group $S_t$ is \emph{max-plus linear}, i.e., 
\begin{equation}\label{stepk}
 S_t[f\vee g]=S_t[f]\vee  S_t[g], \qquad
 S_t[\lambda+ g]=\lambda+ S_t[g]\enspace.
\end{equation}
We shall see that the max-plus basis method exploit these properties to solve the optimal control problem~\eqref{dyuj1}. 
\subsection{Max-plus linear spaces}
A set $\cW$ of functions $\R^d\to \rmax$ is a max-plus linear space 
if for all $\phi_1,\phi_2\in \cW$ and $\lambda\in\R$, the functions
$\phi_1\vee \phi_2$ and $\lambda+\phi_1$ belong to $\cW$. A max-plus
linear space $\cW$ is (conditionally) \emph{complete} if the pointwise
supremum of any family of functions of $\cW$ that is bounded from above
by an element of $\cW$ is finite.

Let $\cB$ be a set of functions $\R^d\to \R$ (\textit{max-plus basis functions}). The complete max-plus (linear) space $\Spanb \cB$ of functions generated by $\cB$ is defined to be the set
of arbitrary linear combinations of elements of $\cB$, in the max-plus
sense, so that an element $\phi$ of $\Spanb\cB$ reads
\(
\sup_{w\in \cB}(a(w) + w) 
\)
for some family $(a(w))_{w\in \cB}$ of elements of $\R_{\max}$.
The (non complete) space $\Span\cB$ is defined in a similar way,
but the linear combination must now involve a finite family, meaning that
$a(w)=-\infty$ for all but finitely many values of $w\in \cB$.
We refer the reader to~\cite{litvinov,ilade,mceneaney-livre} for more background on max-plus linear spaces.

Several choices of basis functions have been considered
in the literature.  
Following~\cite{a5} and~\cite{a6}, we will consider 
here a set $\cB$ consisting of the basis functions of the form
\begin{align}
w_p(x)= -\frac{c}{2}|x|^2+p^T x, \qquad p\in \R^d \enspace ,
\label{e-wp}
\end{align}
where $c$ is a fixed real constant. Hence, an element of $\Spanb \cB$
can be written as 
\begin{equation}\label{keyaz}
\phi(x)=\displaystyle\sup_{p\in \mathbb{R}^d}[-\frac{c}{2}|x|^2+p^T x+a(p)],\quad \forall x \in \R^{d},
\end{equation} 
for some function $a:\R^d\to \rmax$. Recall
that a function $\phi$ is \emph{$c$-semiconvex} if
and only if the function $x\mapsto \phi(x)+\frac{c}{2}|x|^2$ is convex. 
Then, it follows from Legendre-Fenchel duality that the space $\Spanb\cB$ coincides with the space $\cS_c$ of \emph{$c$-semiconvex} (lower semicontinuous)
functions. 

When $c\neq 0$, it is sometimes more convenient to write the 
expansion~\eqref{keyaz} in the form
\[
\phi(x)=\displaystyle\sup_{z\in \mathbb{R}^d}[-\frac{c}{2}|z-x|^2+a'(z)],\quad \forall x \in \R^{d} \enspace.
\]
One can pass from one representation to the other 
by the change of variable $p=cz$.

If $\cW$ is a complete max-plus linear space of functions $\R^d\to\rmax$,
and if $\phi$ is any function $\R^d\to\rmax$, the {\em projection}
of $\phi$ onto $\cW$ is defined to be
\begin{align}
P_{\cW}(\phi):=\max\{\psi \in\cW\mid \psi\leq \phi\} \label{p-primal}
\end{align}
(by writing max, we mean that the supremum element
of the set under consideration belongs to this set,
which follows from the completeness of $\cW$). 

All the previous definitions can be dualized, replacing max by min, and $-\infty$ by $+\infty$. In particular, a complete min-plus linear space is a set
$\cZ$ of functions  $\R^d\to  \R\cup\{+\infty\}$ such that $-\cZ:=\{-w\mid w\in \cZ\}$ is a complete max-plus linear space. Then, we define
the dual projector $P^{\cZ}$ by
\[
P^{\cZ}(\phi):=\min\{\psi \in\cZ\mid \psi\geq \phi\} \enspace,
\]
for all functions $\phi:\R^d\to\R\cup\{+\infty\}$.

\subsection{Max-plus basis methods}\label{subsec-basis}
All these methods approximate the value function at time
$t$ by a finite max-plus linear combination $v_h^t$ of max-plus basis functions, i.e.,
\[ v(x,t)\simeq v_h^{t}(x)=\vee _{i=1}^{q_t}(\lambda_{i}^t + w^{t}_i(x)),
\]
where $\forall i, w^t_i(x) \in \cB$. Typically, $w^t_i(x)= -\frac{c}{2}
|x|^2+p_i^T x$ (see~\cite{a5}), so that the previous approximate
representation is nothing but a discretization of the semiconvex
representation~\eqref{keyaz} of $v(x,t)$.

Then, the coefficients
$\lambda_i^t$ and the functions $w^t_i(x)$ need to be inductively determined.
Let us fix a time discretization
step $\deltat$, such that $T=N\deltat$ for some integer $N$. 
Using the semi-group property, we get
\begin{equation}\label{recur}
 v(\cdot,t+\deltat)=S_{\deltat}[v(\cdot,t)], \qquad t=0,1,\dots, T-\deltat\enspace.
\end{equation}
We require the max-plus basis approximation
$v_h^t$ of $v(\cdot,t)$ to satisfy the analogous relation,
at least approximately:
\begin{eqnarray}\label{formm}
v_h^{t+\deltat}\simeq  S_{\deltat}[v_h^{t}]=
\vee _{i=1}^{q_t}(\lambda_{i}^t + S_{\deltat}[w^{t}_i]),\quad
t=0,\dots,T-\tau
\end{eqnarray}
The various methods differ in the way they address
the two subproblems,
\begin{enumerate}
\item for all $ w^t_i$, replace $S_{\deltat}[w^t_i]$ 
by an easily computable (accurate enough) approximation
$\tilde{S}_{\deltat}(w^t_i)$;
\item Project each $\tilde{S}_{\deltat}(w^t_i)$, or
$\vee _{i=1}^{q_t}(\lambda_{i}^t + S_{\deltat}[w^{t}_i])$,
to the space of finite max-plus linear combinations of basis functions. 
\end{enumerate}
The first subproblem is the simplest one: 
computing $S_{\deltat}[w^t_i]$ is equivalent to solving
an optimal control problem, but the horizon $\deltat$ is small,
and the terminal reward $w^t_i$ (typically a quadratic function) has
a regularizing and a ``concavifying'' effect, which implies
that the global optimum can be accurately approached (by reduction to a convex
programming problem), leading to various approximations
with a consistency error of $O(\tau^r)$, with $r=3/2, 2$, or sometimes better,
depending on the scheme, see~\cite{mceneaney-livre,a6,theseasma}.

The second step is the critical one,
in particular, the accuracy of the method is
limited by the \emph{projection error}~\cite{a6},
i.e., the maximal distance between every function
$S_{\deltat}[v_h^{t}]$ or $ S_{\deltat}[w^{t}_i]$ and its best approximation
by a max-plus linear combination of basis functions. In a number
of methods, including the original one~\cite{a5}, the approximation
$v_h^t$ is such that, assuming
that the approximation $\tilde{S}_{\deltat}(w_i^t)$ of $S_{\deltat}(w_i^t)$ is
exact,
\begin{align}
v_h^t\leq v(\cdot,t)\label{e-islower}
\end{align}
so that
\(
v_h^t \leq P_{\cW_t}v(\cdot,t)
\)
where $\cW_t=\Span\{w_i^t\mid  1\leq i\leq q_t\}$, and $P_{\cW_t}$
is defined as in~\eqref{p-primal}. Then, the approximation error 
in the $L_p$ norm, $\epsilon_p:= \|v(\cdot,t)-v_h^t\|_p$ satisfies 
\begin{align}
\epsilon_p \geq \|v(\cdot,t)-P_{\cW_t}v(\cdot,t)\|_p
\label{e-epsilon}
\end{align}

Similarly, in the max-plus finite element method of Akian, Gaubert and Lakhoua~\cite{a6}, the approximation $v_h^t$ is computed recursively by
\[
v_h^t= P^{\cZ_t}P_{\cW_t} \tilde{S}_\deltat(v_h^{t-\tau}) \enspace, 
\]
where at each step, we also have a (dual) min-plus linear space $\cZ_t$ generated by finitely many test-functions. Then, it follows from~\cite{a6} that the sup-norm projection error cannot be expected to be smaller than
\[
\| v(\cdot,t-\tau)
- P_{\cW_t}v(\cdot,t-\tau)\|_\infty
+
\| v(\cdot,t-\tau)
- P^{\cZ_t}v(\cdot,t-\tau)\|_\infty \enspace .
\]
Moreover, it is shown in~\cite{a6} that this estimate is tight, meaning
that the total error of the method, $\|v(\cdot,T)-v_h^T\|_\infty$
is of order at most $N$ times the previous sum,
up to a term depending on the quality of approximation 
of $S_\tau[w]$ by $\tilde{S}_\tau[w]$.

Finally, in the curse of dimensionality method of McEneaney~\cite{curseofdim}, the basis functions are quadratic forms
and the semi-group $S_{\deltat}$ is approximated by the pointwise
supremum of semi-groups associated to linear-quadratic control problems
\begin{equation}\label{whhh}
\tilde S_{\deltat}[\phi]=\vee_{m=1}^M S_{\deltat}^m[\phi] \enspace.
\end{equation}
The number of basis functions of the linear combination grows by a factor of $M$ at each step. Then, the accuracy of the method is still limited by the projection error, which arises when pruning the less useful basis functions.

\section{Curse of dimensionality for semiconvex based approximations}
\label{sec-curseunavoidable}

As discussed above, the main source of inaccuracy of max-plus basis
methods is the projection error~\eqref{e-epsilon}. In this section, we 
give an asymptotic estimate of the optimal projection error
as the number of basis functions tends to infinity, in the special
case in which the space $\cW_t$ consists of quadratic functions~\eqref{e-wp},
as in the max-plus basis method~\cite{a5}, or in the $P_2$ type
finite element method of~\cite{a6,theseasma}.

Let $c\in \R$, $\epsilon>0$ and $\psi(x):\R^d \rightarrow \R$ be a $(c-\epsilon)$-semiconvex function. It is approximated
by a given number $n$ of semiconvex basis functions
$-\frac{c}{2}|x|^2+p^T x$:
\begin{align}
\psi(x)\simeq \displaystyle \tilde{\psi}(x):=\max_{i=1,\dots,n} \{\
-\frac{c}{2}|x|^2+p_i^T x
+a(p_i)\} 
\enspace .\label{e-approx}
\end{align}
We are interested in the 
$L_1$ or $L_\infty$ approximation error
\[
\epsilon_1:=\int_X |\psi(x)-\tilde{\psi}(x)|dx,\quad\text{or}\quad
\epsilon_\infty:=\sup_{x\in X} |\psi(x)-\tilde{\psi}(x)|
\]
for some suitable full dimensional compact convex subset $X\subset\R^d$.
We shall compute these errors when $\psi=v(\cdot,t)$, and so,
for reasons discussed in Section~\ref{subsec-basis} (see~\eqref{e-islower}), we shall require
that $\tilde{\psi}\leq \psi$.
We denote by $\delta_{X,n}^1(\psi,c)$ (resp.\ $\delta_{X,n}^{\infty}(\psi,c)$) the minimal $L_1$ (resp.\ $L_\infty$) approximation error on $X$ of $\psi(x)$ by $\tilde{\psi}(x)$ as in~\eqref{e-approx}, by $\psi_x''$ the hessian matrix of $\psi$ at point $x$, and by $I_d$ the identity matrix of size $d$.

The next two theorems imply
that whatever computation scheme is chosen for the coefficients
$a(z_i)$, the approximation error is necessarily subject to a curse
of dimensionality. 
\def\localskip{\\[1mm]}
\begin{theo}[$L_1$ approximation error]\label{thh3}
Let $c\in \R$, $\epsilon>0$ and let $X\subset\R^d$ denote any full dimensional compact convex subset. If $\psi(x):\R^d \rightarrow \R$ is $(c-\epsilon)$-semiconvex of class $\mathcal{C}^2$, then, there exists a constant $\alpha_1>0$ depending only on $d$ such that
\begin{equation*}\label{L1}
 \delta_{X,n}^1 (\psi,c) \sim \frac{\alpha_1}{n^{\frac{2}{d}}}\Big(\int_X (\det (\psi_x''+c I_d))^{\frac{1}{d+2}} \, dx\Big)^{\frac{d+2}{d}}
\text{ as }n\to\infty
\enspace .\localskip
\end{equation*}
\end{theo}
\begin{theo}[$L_\infty$ approximation error]\label{thh4}
Let $c$, $\epsilon$, $X$ and $\psi(x)$ be as in Theorem~\ref{thh3}. Then
there exists a constant $\alpha_2>0$ depending only on $d$
such that
\begin{equation*}
 \delta_{X,n}^{\infty} (\psi,c) \sim \frac{\alpha_2}{n^{\frac{2}{d}}}\Big( \int_X (\det (\psi_x''+c I_d))^{\frac{1}{2}} dx \Big)^{\frac{2}{d}}
\text{ as }n\to\infty
\enspace .\localskip
\end{equation*}
\end{theo}
The proof of these theorems is sketched in the next section, it builds on analogous methods and results of Gruber~\cite{a2},~\cite{b2}, concerning 
the approximation of smooth strictly convex bodies by circumscribed
polytopes, the constants $\alpha_1$ and $\alpha_2$,
which grow slowly with $d$, already appeared there.

The following theorem is a direct corollary:
\begin{theo}\label{theo-unavoidable}
Assume that $\psi:=v(\cdot,T)$ is the value function, and that it is $\mathcal{C}^2$ and $c$-semiconvex. Then, for any max-plus basis method providing an approximation from below of of the value function by a supremum
of $n$ quadratic functions (see~\eqref{e-approx}), the
$L^1$ error $\epsilon_1$ and $L^\infty$ error $\epsilon_\infty$ of approximation of the value function are both
$\Omega\big( \frac{1}{n^{2/d}} \big)$
as $n\to\infty$.
\end{theo}

Besides, the estimates of Theorems~\ref{thh3} and~\ref{thh4}
confirm that when $n$ is sufficiently large, it is more interesting to choose the smallest $c$ such that $\psi(x)+\frac{c}{2}|x|^2$ is convex. Note
also that the integral term can be small if the Hessian of $\psi$
is nearly constant and close to $-cI_d$ (attenuation of the curse of dimensionality).

\section{Sketch of proof of Theorems~\ref{thh3} and~\ref{thh4}}\label{sec-proof}
In this section, we sketch the proof of Theorems~\ref{thh3} and~\ref{thh4}. Let $c$, $\epsilon$, $X$ and $\psi(x)$ satisfy the conditions of these theorems. Note that approximating $\psi$ as in~\eqref{e-approx} is equivalent to approximating the strongly convex function $\psi(x)+\frac{c}{2}|x|^2$ by $n$ of its affine minorants.
Hence, it suffices to prove Theorems~\ref{thh3} and~\ref{thh4} when $c=0$,
$\psi$ is strongly convex, and is approximated by the supremum of $n$ of its affine minorants. We denote by $\nabla \psi$ the gradient operator of $\psi$, by $\hat x$ the point $(x,\psi(x))\in \mathbb{R}^{d+1}$ and by $\psi_x''(\cdot)$ the quadratic form determined by the hessian matrix $\psi_x''$.
Let us first recall some results on the approximation of convex bodies by circumscribed polytopes.

Gruber proved in~\cite{a2},~\cite{b2} that for any convex body $C\subset \mathbb{R}^{d+1}$ with $\mathcal{C}^2$ boundary and positive Gaussian curvature $\kappa_C>0$, there are two constants $\alpha_1$ and $\alpha_2$ depending only on $d$ such that as $n\rightarrow \infty$:
\begin{eqnarray*}
&&\displaystyle \delta_n^V(C)\sim \alpha_1\Big(\int _{\partial \mathcal{C}} \kappa_C(x)^{\frac{1}{d+2}} d \sigma (x) \Big)^{\frac{d+2}{d}} \frac{1}{n^{\frac{2}{d}}},
\\
&&\displaystyle \delta_n^H(C)\sim \alpha_2\Big(\int_{\partial \mathcal{C}} (\kappa_C(x))^{\frac{1}{2}}d \sigma (x)\Big)^{\frac{2}{d}} \frac{1}{n^{\frac{2}{d}}}\enspace.
\end{eqnarray*}

Here $\delta_n^H(C)$ and $\delta_n^V(C)$ are respectively the minimal distance with respect to Hausdorff and $L_1$ metric between $C$ and any circumscribed polytope with $n$ facets, $\sigma$ is the ordinary surface area measure on $\partial C$. 
Moreover, we have the following asymptotic estimates for the constants $\alpha_1$ and $\alpha_2$:
\[%
 \alpha_1\sim \frac{d+1}{\pi e}, \alpha_2\sim \frac{1}{2\pi}\big (\Gamma (\frac{d}{2}+1)\vartheta_d\big)^{\frac{2}{d}},\quad\mathrm{as}\quad d\rightarrow \infty,
\]%
where $\vartheta_d$ is estimated as~\cite{b1}:
\[%
\deltat_d \leq \vartheta_d \leq d \log d+d \log(\log d)+5d,
\;\mathrm{with}\quad\deltat_d \sim \frac{d}{e\sqrt e}\enspace.
\]%

Both of the proofs partition $\partial C$ into finitely many pieces associated
to a family of supporting planes. For each supporting plane, there is a corresponding strongly convex function whose graph is a piece of $\partial C$. Then,
the volume of the difference between $C$ and a circumscribed polytope
can be estimated by computing the $L_1$ norm of the difference of
this strongly convex functions with some of its piecewise affine 
lower bounds. For the Hausdorff metric case, some results regarding
the optimal covering of a manifold by geodesic discs are used. 
We next apply the same techniques to our problem.

First of all, we recall the definition of the~\textit{Bregman distance}.

\begin{de}[{\cite{MR0215617}}]
For any two points $x$ and $y$ of $X$, the \textit{Bregman distance} $D_\psi(\cdot;\cdot)$ from $x$ to $y$, associated to a strongly convex and differentiable function $\psi$, is defined by
\begin{equation}
 D_\psi(x;y)=\psi(x)-\psi(y)-\nabla \psi(y)^T (x-y)\enspace .
\end{equation}
\end{de}
The Bregman distance is positive definite ($D_\psi(x,y)\geq 0$ and the equality
holds if and only if $x=y$), but it may not be symmetric.

The proofs follow essentially the same steps as for convex bodies. We give just an outline of the proof without going into details. 

To prove Theorem~\ref{thh3}, we need an asymptotic formula 
for optimal quantization: 
\begin{theo}[{\cite{a2}}]\label{zador} 
 Let $J\subseteq \mathbb{R}^d$ be Jordan measurable with positive volume $v(J)>0$, and $q$ a positive definite quadratic form on $\mathbb{R}^d$. Then as $m \rightarrow \infty$:
$$\inf_{S\subseteq \mathbb{R}^d, |S|=m} \int _J \min_{t\in S} \{q(s-t)\} ~ds \sim 2\alpha_1 v(J)^{\frac{d+2}{d}} (\det q)^{\frac{1}{d}} \frac{1}{m^{\frac{2}{d}}}\enspace. $$ 
\end{theo}

\textit{Sketch of proof of Theorem~\ref{thh3} for $c=0$}

Let $\lambda>1$, for each $p\in X$, there is an open convex neighborhood $U\subset X$ such that:
$$
\frac{1}{\lambda}\psi_u''(x)\leq \psi_p''(x)\leq \lambda \psi_u''(x), \quad\forall u\in U,\quad \forall x\in \R^d\enspace.
$$
Then for every $x, y\in U$, the Bregman distance $D_\psi(x;y)$ is bounded below and above as:
\begin{equation}\label{equ12}
\frac{1}{2\lambda} \psi_p''(x-y)\leq D_\psi(x;y)\leq \frac{\lambda}{2} \psi_p''(x-y) \enspace.
\end{equation}

We choose finitely many points $p_1,p_2,\dots,p_m$ with respective
neighborhoods $U_1,\dots,U_m$ covering the compact set $X$. 
One may then dissect the integral
$\int_{X} [\min_{y\in S} D_\psi(x;y)] dx$
on smaller pieces $\{J_i\subset U_i, i=1\dots,m\}$ with $\{J_i,i=1,\dots,m\}$ Jordan measurable. Using the asymptotic formula of Theorem \ref{zador} and some other arithmetic inequalities as in~\cite{b2}, the theorem can be deduced.

\textit{Sketch of proof of Theorem~\ref{thh4} for $c=0$}

Let $\hat X$ be the graph of $\psi(x)$ on $X$. $\hat X$ is a $d$-dimensional (Riemannian) manifold of class $\mathcal{C}^2$ with metric of class $\mathcal{C}^0$. For each $\hat x, \hat y \in \hat X$, one may define the 
Riemannian metric between $\hat x$ and $\hat y$, $\gamma(\hat x,\hat y)$,
by: 
\[%
 \gamma(\hat x,\hat y)=\inf \{\int_0^1 \psi''_{u(t)}(\dot{u}(t))^{\frac{1}{2}} \,dt|u(t)\in C^1\!, u(0)=x,u(1)=y\}.
\]%

To prove Theorem~\ref{thh4}, we need an asymptotic formula on the minimum covering. The next lemma is a special case of Lemma~1 in~\cite{a2}:
\begin{lemma}[Compare with Lemma 1 in~{\cite{a2}}]\label{lmcom}
For $\rho >0$, let $n(\hat X,\rho)$ be the minimum number of discs of radius $\rho$ with respect to the Riemannian metric needed to cover $\hat X$. Then:
\begin{equation}\label{haus}
n(\rho)\sim \alpha_2^{d/2}\Big (\int _X (\det \psi_x'')^{\frac{1}{2}} dx\Big )\rho^d, \quad\mathrm{as }\quad \rho \rightarrow 0\enspace. \end{equation}
\end{lemma}
Given a similar result of minimum covering under Euclidean metric of Hlawka~\cite{MR0030547}, this lemma essentially proves the equivalence between the Riemannian metric on $\hat X$ and the Euclidean metric on $X$. Since our problem is to minimize the maximum radius of Bregman balls covering the manifold $\hat X$, one last thing to be proved is the equivalence between the Riemannian distance and the Bregman distance. Indeed, we prove that:
\begin{equation}\label{equ14}
 \exists M>0, \forall x,y \in X, 
\frac{\gamma(\hat x, \hat y)}{D_\psi(x;y)}\leq M\enspace,~~~~~~~~~~~~~~~~~~
\end{equation}
\begin{equation}\label{equ15}
\begin{array}{ll}
\forall x,y \in U_l, &\gamma(\hat x, \hat y)< \dist_\gamma (\hat x,\partial \hat U_l), \Rightarrow \\
&\frac{1}{2\lambda^4}\gamma (\hat x,\hat y)^2 \leq D_\psi(x,y)\leq \frac{\lambda^4}{2} \gamma (\hat x,\hat y)^2\enspace.
\end{array}
\end{equation}
Here $\dist_\gamma(\hat x, \partial \hat U_l)=\min\{\gamma(\hat x,\hat y):y\in \partial \hat U_l\}$.
Using the minimum covering asymptotic estimation \eqref{haus}, the equivalence between Bregman distance and Riemannian distance \eqref{equ14} and \eqref{equ15}, the desired theorem is established as in~\cite{a2}.

\section{Primal curse of dimensionality free method}\label{sec3}
We consider the optimal control problem for switched linear system studied in~\cite{curseofdim} (see also~\cite{a7,mccomplex}). Let ${\cal{M}}=\{1,2,\dots,M\}$.
\begin{equation}
 V(x)=\sup_{\omega\in \setomega}\sup_{\mu \in \cal{D}_{\infty}} \sup_{T<\infty} \int_{0}^T L^{\mu_t}(\xi_t)-\frac{\gamma^2}{2}|\omega_t|^2 dt, 
\end{equation}
where
 \begin{eqnarray}
  &&L^{\mu_t}(x)=\frac{1}{2}x^TD^{\mu_t}x+(l_1^{\mu_t})^Tx+\alpha^{\mu_t}, \\
&&{\cal{D}}_{\infty}=\{\mu: [ 0,\infty )  \rightarrow \cal{M}: \mathrm{measurable} \},\\
&&{\setomega} \doteq {L}_{2}^{\mathrm{loc}}([\mathrm{0},\infty);\mathbb{R}^k),
 \end{eqnarray}
and the state dynamics are given by
\begin{equation}
 \dot \xi=A^{\mu_t}\xi+l_2^{\mu_t}+\sigma^{\mu_t}\omega_t,\xi_0=x,
\end{equation}
where $\sigma^m$ and $\gamma$ are such that $\Sigma^m=\frac{1}{\gamma^2}\sigma^m(\sigma^m)^T$, $\forall m\in {\cal{M}}$.
The corresponding HJB PDE is:
\begin{equation}\label{hjb}
 0=-H=-\max_{m\in {\cal{M}}} \{H^m(x,\nabla V)\},
\end{equation}
where $H^m$ has the form:
\begin{equation}\label{Hm}
\!\!\!\!\!\begin{array}{rcl}
 H^m(x,p)&\!=\!&\frac{1}{2}x^T D^m x+\frac{1}{2}p^T \Sigma^m  p +\\
&&\quad (A^m x)^Tp+(l_1^m)^Tx+(l_2^m)^Tp+\alpha^m\enspace.
\end{array}
\end{equation}
Under certain technical assumptions~\cite{curseofdim} which we will not repeat here, the function $V$ is finite,
it is a viscosity solution of~\eqref{hjb}, and it is given by
$
 V=\lim_{T\rightarrow \infty} S_T[0]
$
where $(S_t)$ is the Lax-Oleinik semi-group of the Hamilton-Jacobi equation~\eqref{guih} for $H$ defined in~\eqref{hjb}.
In~\cite{curseofdim}, the value function is approximated by a max-plus sum of quadratic functions, and the approximated semi-group is propagated in a dual space.  We next introduce a variant, which we call the~\textit{primal} curse of dimensionality free method: it is equivalent if no trimming is performed, but it avoids the use of dual representations.

\subsection{Approximate propagation}%
The basis functions are allowed to be all of the quadratic functions smaller than the value function.  
Let $S_{\deltat}^m$ be the evolution semi-group of the Hamilton-Jacobi equation~\eqref{guih} for $H_m$ defined in~\eqref{Hm}, $m\in {\cal{M}}$.
We approximate $S_{\deltat}[\phi]$ by:
\begin{equation}
 S_{\deltat}[\phi]\simeq \tilde{S}_{\deltat}[\phi]=\vee_{m\in {\cal{M}}} S_{\deltat}^m[\phi],
\end{equation}
which for $\phi$ quadratic still yields a maxima of quadratics.

\subsection{Computation of a single semi-group operator}
The propagation of a quadratic function $\phi$ by $S_{\deltat}^m$ reduces
to solving a differential Riccati equation (DRE).
Moreover, it is well-known that one can recover the solution
of a DRE from a system of
Hamiltonian linear differential equations (see, e.g.,~\cite{MR0357936}).
Suppose there are only quadratic terms, i.e., $l_1^m=0,l_2^m=0, \alpha^m=0$. Let $\phi(x)=\frac{1}{2}x^TP_0x$, then $S_{t}^m[\phi](x)=\frac{1}{2}x^T P_t x$, where $P_t=Y_t{X_t}^{-1}$ and $(X_t,Y_t)$ are the solution of:
\begin{equation}\label{linearsystem}
\left\{
\begin{array}{ll}
\left(
\begin{array}{c}
\dot X   \\
\dot Y
\end{array}
\right)
=
\left(
\begin{array}{ll}
-A^m & -\Sigma ^m\\
D^m & (A^m)^T
\end{array}
\right )
\left (
\begin{array}{c}
X \\
Y
\end{array}
\right)\\ \\
\begin{array}{c}
 X(0)=I_d, \quad Y(0)=P_0  
\end{array}
\end{array}
\right. \enspace.
\end{equation}
We denote by $\mathcal{A}$ the matrix coefficient in the above linear system.
Note that the invertibility of $X(t)$ can be derived from the fact that the value function $ V$ is finite. Given a fixed time step $\deltat>0$, the fundamental solution $\exp(\mathcal{A}\deltat)$ of the previous linear system satisfies:
\[%
 \left(\begin{array}{l}
  X_t\\Y_t
 \end{array}
\right)=\exp (\mathcal{A}\deltat) \left(\begin{array}{l}
  I_d\\P_0
 \end{array}
\right)\enspace.
\]%
In the presence of linear or constant terms in the control system or in the quadratic function, the problem can be easily transformed into a purely quadratic one by adding a constant state variable. The above analysis shows that, given a fixed propagation time $\deltat$, computing $S^m_{\deltat}[\phi]$, for
every quadratic form $\phi$, reduces
to a matrix multiplication and an inverse operation, which can be done
in $O(d^3)$ incremental time.

\subsection{Propagation and curse of complexity}

We choose a time discretization step $\deltat$, a number of steps $K$ and an initial function $\phi_{0}(x)$. Let $\phi_l(x)=\vee_{j \in J} \phi_l^j(x)$ be the approximation at step $l$. The iteration formula is given by:
\[%
 \phi_{l+1}(\cdot)=\vee_{m\in {\cal{M}}} S_{\tau}^m[\phi_{l}(\cdot)]=\vee_{m\in {\cal{M}},j\in J} S_{\deltat}^m[\phi_l^j(\cdot)]\enspace.
\]%
The computational growth in the space dimension is cubic, as shown in the above subsection. However, the number of quadratic forms grows by a factor of $M$ at each iteration. This \emph{curse-of-complexity} issue also occurs for the dual method in~\cite{curseofdim}. Some SDP relaxation based pruning method was proposed in~\cite{a7} to reduce the number of quadratic forms. We next discuss improvements of this pruning methods, still partly SDP based, but now exploiting 
the combinatorial nature of the problem.

\section{Reduction of pruning to $k$-center and $k$-median problems for a Bregman type distance}\label{sec4}
\subsection{Formulation of the pruning problem}
We first give a general formulation for the pruning problem appearing in max-plus basis methods. Let $F=\{1,2,\dots,n_f\}$.
Let $\phi(x): \mathbb{R}^d\rightarrow \mathbb{R}$ be a max-plus sum of $n_f$ basis functions: $\displaystyle \phi(x)=\vee_{j\in F} \phi_j(x)$. Let $0<k<n_f$ be a fixed integer. 
The problem is to approximate $\phi(x)$ by keeping only $k$ basis functions. To measure the approximation error, we introduce a Bregman type distance $\dist_{\phi}(x;j)$ between each point $x\in \R^d$ and 
each basis function $\phi_j(\cdot)$, such that:
\begin{eqnarray*}
\!\!\!\!\!\!\!\!\!\!&&\forall x\in \R^d,~ \exists j_0\in F,~\mathrm{s.t.}~ \dist_{\phi}(x;j_0)=0\enspace;\\
\!\!\!\!\!\!\!\!\!\!&&\forall x\in \R^d, i, j\in F,~\dist_{\phi}(x;i)\leq \dist_{\phi}(x;j) \Leftrightarrow \phi_j(x)\leq \phi_i(x) \,.
\end{eqnarray*}
In other words, the distance $\dist_{\phi}(x;j)$ measures the loss at point $x$ caused when approximating $\phi(\cdot)$ by $\phi_j(\cdot)$.  
For example, the simplest choice is to let $\dist_\phi(x;j)=\phi(x)-\phi_j(x)$. Consider a compact set $X\subset \R^d$ on which we measure the loss. One may minimize the total loss ($L_1$ metric) or the maximal loss ($L_\infty$ metric) on $X$.

\subsubsection{$L_1$ metric and $k$-median problem}
\begin{equation}\label{zfter3}
 \delta_{k}^1(\phi)=\min_{\substack{S\subset F,|S|=k}} \int_{X} [\min_{j\in S} \dist_\phi(x;j)] dx\enspace.
\end{equation}
\subsubsection{$L_\infty$ metric and $k$-center problem}
\begin{equation}\label{forlu2b}
 \delta_{k}^{\infty}(\phi)=\min_{\substack{S\subset F,|S|=k}} \max_{x\in X} [\min_{j\in S} \dist_\phi(x;j)]\enspace.
\end{equation}

We recognize in \eqref{zfter3} and \eqref{forlu2b} the classical $k$-median and the $k$-center facility location problem with continuous demand area and discrete service points. The facility location problem, discrete or continuous, is known to be $NP$-hard even with euclidean distance. Besides, we
remark that a subproblem of Problem~\eqref{zfter3} is the
volume computation for polytopes, which is known to be $\#P$-hard. 
To the best of our knowledge, the only few references that discuss this general class of location problem replace the continuous demand with a discrete one with large number of points, see~\cite{MRDrezner}.
In the following, we consider a specific case and propose a method based on SDP relaxation to generate discrete points.

\subsection{Pruning methods}\label{pruningmethod} 
We assume that all basis functions are quadratic: $\phi_j(x)=\frac{1}{2}x^TA_jx+b_j^Tx+\frac{1}{2}c_j, \forall j\in F$. We normalize the distance function as in~\cite{a7}, i.e.,
$\dist_{\phi}(x;j)=(\phi(x)-\phi_j(x))/(1+|x|^2)\enspace.$
\subsubsection{'Sort upper bound'~\cite{a7}}
The first method was introduced in~\cite{a7}. Roughly speaking, to each basis function $\phi_j(x)$ we associate an \textit{importance metric} :
\begin{equation}\label{impmetric}
 \nu_j=\max_{x\in \mathbb{R}^d} \min_{j'\neq j}(\phi_j(x)- \phi_{j'}(x))/(1+|x|^2)\enspace.
\end{equation}
Then $\nu_j$ is the normalized $L_\infty$ error caused by pruning the function $\phi_j(x)$.
In some sense the bigger $\nu_j$ is, the more useful the function $\phi_j(x)$ is. In particular, when $\nu_j\leq 0$ the function $\phi_j(x)$ is dominated by the others and it can be pruned without generating any approximation error.
Let
\begin{displaymath}
 Q_j^{j'}=\frac{1}{2}\big[\begin{array}{ll}
        c_j-c_{j'} & b_j^T-b_{j'}^T \\
b_j-b_{j'} & A_j-A_{j'}
       \end{array}\big]=Q_j-Q_{j'}\enspace.
\end{displaymath}
The problem \eqref{impmetric} is equivalent to:
\begin{equation}\label{equi}
\! \nu_j\!\!=\!\!\!\!\max_{\substack{\nu \in \R ; y \in \R^{d+1}}}\!\{ \nu: y_1\neq 0; \|y\|=1;y^T Q_j^{j'} y \geq \nu, \forall j'\neq j\}.
\end{equation}
This nonconvex QCQP(quadratically constrained quadratic program)~\cite{MR2061575} has
its SDP relaxation given by:
\begin{equation}\label{sdprel}
\!\overline {\nu}_j\!\!=\!\!\!\!\max_{\substack{z\in \R,Y\succeq 0\\ y \in \R^{d+1}}}
\!\left\{\nu \left |
\begin{array}{ll}
 Y_{11}>0; \quad \Tr(Y)=1;\quad Y \succeq yy^T;\\
\Tr(Y  Q_i^j)\geq \nu, \quad \forall j\neq i .\\
\end{array}\right\}\right..
\end{equation}
Then $\overline \nu_j$ is an upper bound of the importance metric $\nu_j$. Finally the \textit{sort upper bound} method consist in sorting all the upper bounds $\{\overline \nu_j,j\in F\}$ and picking up the $k$ first ones.

\subsubsection{'Sort lower bound'}

The SDP relaxation~\eqref{sdprel} provides not only an upper bound on the importance metric but also a rather simple way to generate feasible solutions. 

Suppose $(\overline Y,\overline y)$ is a solution of program \eqref{sdprel}. We use the randomization technique~\cite{MR1801497} to get feasible points: we pick $y$ as a Gaussian random variable with $y\sim \mathcal{N}(\overline y, \overline Y-\overline y \overline y^T)$. 
Then over this distribution, the constraints in \eqref{equi} are satisfied on average. 
By sampling $y$ a sufficient number of times, we get a $y$ such that the inequality constraints in \eqref{equi} are all satisfied. Then, setting $x=(y_2/y_1,\dots,y_{d+1}/y_1)^T$ provides a lower bound of \eqref{impmetric}.
The proposed procedure provides in practice a good lower bound,
although there is no theoretical guarantee in the present generality.

Then, the \textit{sort lower bound} method proceeds as follows.
Fixing an integer $N>0$, for each basis function $\phi_j$, we resolve the SDP program \eqref{sdprel}, get $N$ feasible points $x\in \R^d$ using the above randomization technique and put them into a set $X'$. At the end we get a discrete set $X'$ and for each basis function $\phi_j(x)$ we calculate its lower bound ${\underline \nu}_j$ by:
$${\underline\nu}_j=\max_{x\in X'} \min_{j'\neq j}(\phi_j(x)- \phi_{j'}(x))/(1+|x|^2)\enspace.$$ 
Finally we sort all of the lower bounds $\{\underline{\nu}_j,j\in F\}$ and keep the $k$ first ones.

Following the above randomization technique, we get a discrete set $X'$ which in some sense reflect rather well the importance of each basis function.
We replace the compact set $X$ by this discrete set $X'$ and seek to minimize the total loss on $X'$. This gives the discrete $k$-median problem:
\begin{equation}\label{diskm}
\delta=\min_{\substack{S\subset F,|S|=k}} \sum_{x\in X'} [\min_{j\in S} \dist_{\phi}(x;j)]\enspace.
\end{equation} 
This central problem in combinatorial optimization has seen a succession of papers designing approximations algorithms. Our two last pruning methods are merely two heuristics for the $k$-median problem~\eqref{diskm}.

\subsubsection{'J-V facility location'}
Lin and Vitter~\cite{Lin:1992:EMP:129712.129787} proved that the constant factor approximation for general $k$-median problem is $NP$-hard. For metric distance, Jain and Vazirani~\cite{MR1851303} proposed a primal-dual 6-approximation algorithm. This algorithm is interesting not only due to its constant factor,
but also because it is combinatorial (there is no need to solve a linear program).  

\subsubsection{'greedy facility location'}
The fourth method is the greedy heuristic. Remember that the function
to be minimized in the facility location problem is supermodular,
which implies that the greedy heuristic has a bound estimate
(even without the triangular inequality on the distance function).
Let $\delta_G$ be the value of a particular solution constructed by the greedy heuristic, then we have~\cite{MR0503866}:
\(%
\delta_G\leq (1-\alpha^k)\delta+\alpha^k(\max_{j\in F} \sum_{x \in X'} \dist_{\phi}(x;j)),
\)%
where $\alpha=\frac{k-1}{k}$.
The execution time of the greedy heuristic is $O(km)$.

\section{Experimental results}\label{sec-exp}
\subsection{Problem instance}
To compare with the dual max-plus basis method, 
we use the instance of~\cite{a7} originating
from H-infinity control, in which the parameters where
chosen to exhibit a complex behavior.
The state dimension and the switch number are both 6. 
The overpruning threshold is also the same: we keep $k(i)=20+6i$ basis functions at step $i$.

Without the exact value function, we do not have a direct error estimation. Recall that the value function $V$ is the unique viscosity solution of the following HJB equation:
\begin{equation}\label{hjb2}
 0=-H=-\max_{m\in {\cal{M}}} \{H^m(x,\nabla V)\},
\end{equation}
where $H^m$ is defined in~\eqref{Hm}. The value of Hamiltonian is then used to measure the approximation.
\subsection{Numerical results}
 All of our results\footnote{The code was mostly written in Matlab (version 7.11.0.584), calling YALMIP (version 3) and SeDuMi (version 1.3) for the resolution of SDP programs. 
The computation of the distance function and Jain \& Vazirani's primal dual algorithm were written in C++. 
The results were obtained on a single core of an Intel quad core running
at 2.66GHz, with 8Gb of memory. }
 are shown along the $x_1$-$x_2$ axes with the 4 other coordinates of $x$ set to $0$.

\begin{figure}[thpb]
      \centering
      \includegraphics[scale=0.4]{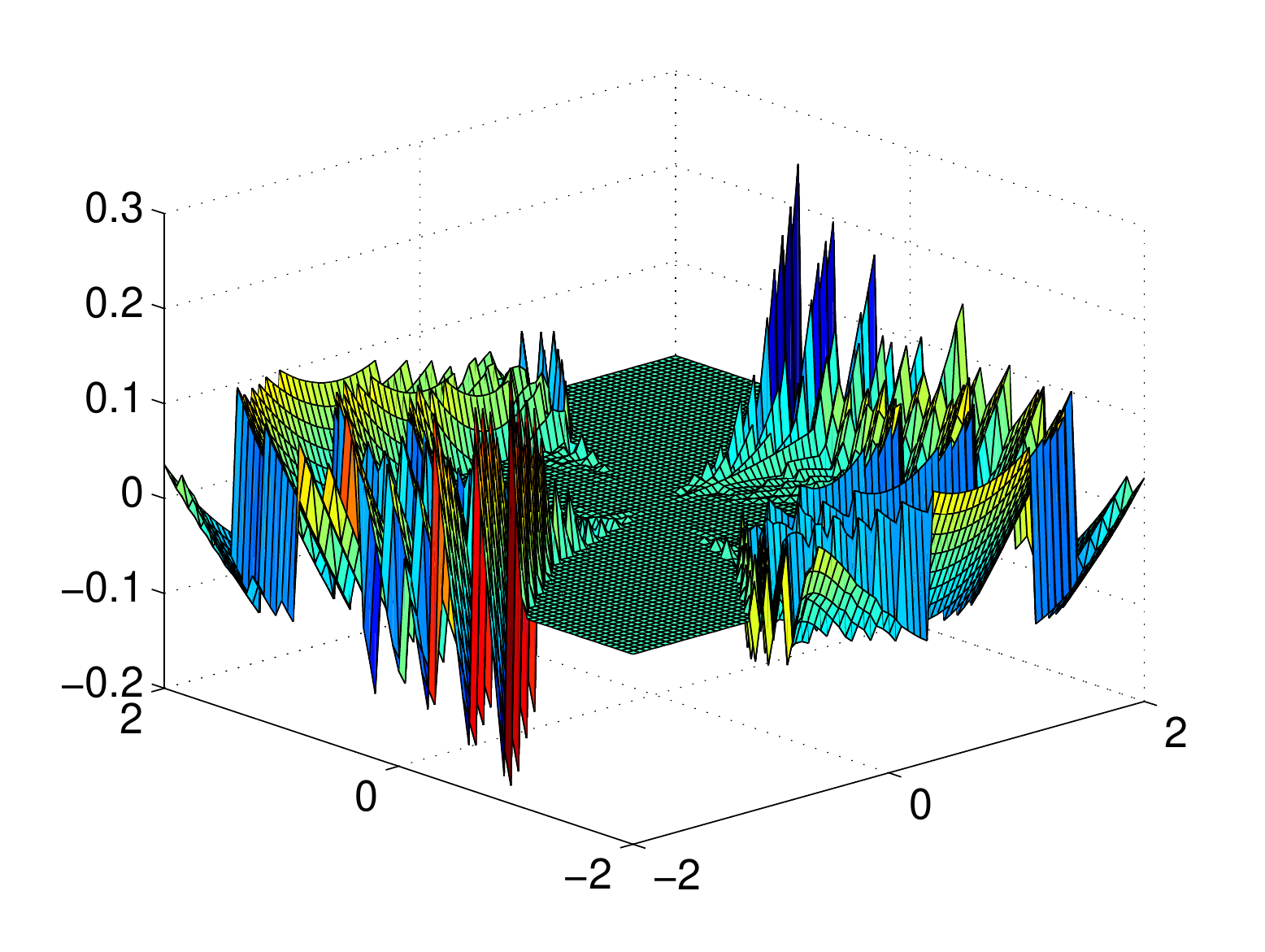}
      \caption{Backsubstitution error on the $x_1$,$x_2$ plane, $\deltat=0.2$, with \textit{sort upper bound} pruning method}
      \label{fl1}
   \end{figure}

For comparison with~\cite{a7}, we first take the same time-discretization step-size $\deltat=0.2$, the same iteration steps 25 and the same pruning method
(\textit{sort upper bound}). %
Figure~\ref{fl1} shows the value of Hamiltonian $H$ at the end of 25 iterations. Comparing with the error plot shown in~\cite{a7}, which is in the same scale but has a peak of error of order $1$ (versus $0.3$ here), we see that
the primal max-plus basis method yields a small improvement.

\begin{figure}[thpb]
      \centering
      \includegraphics[scale=0.4]{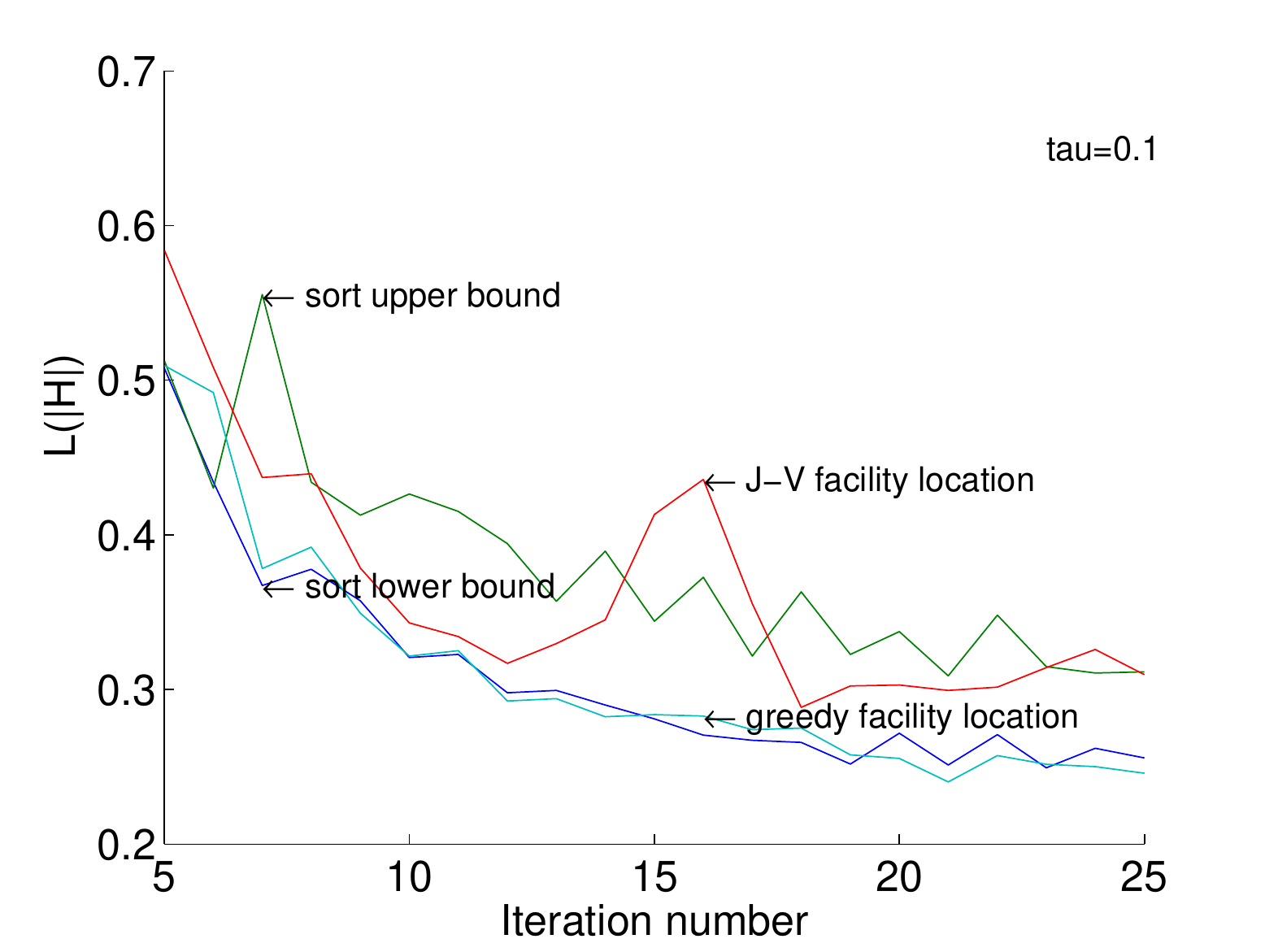}
      \includegraphics[scale=0.4]{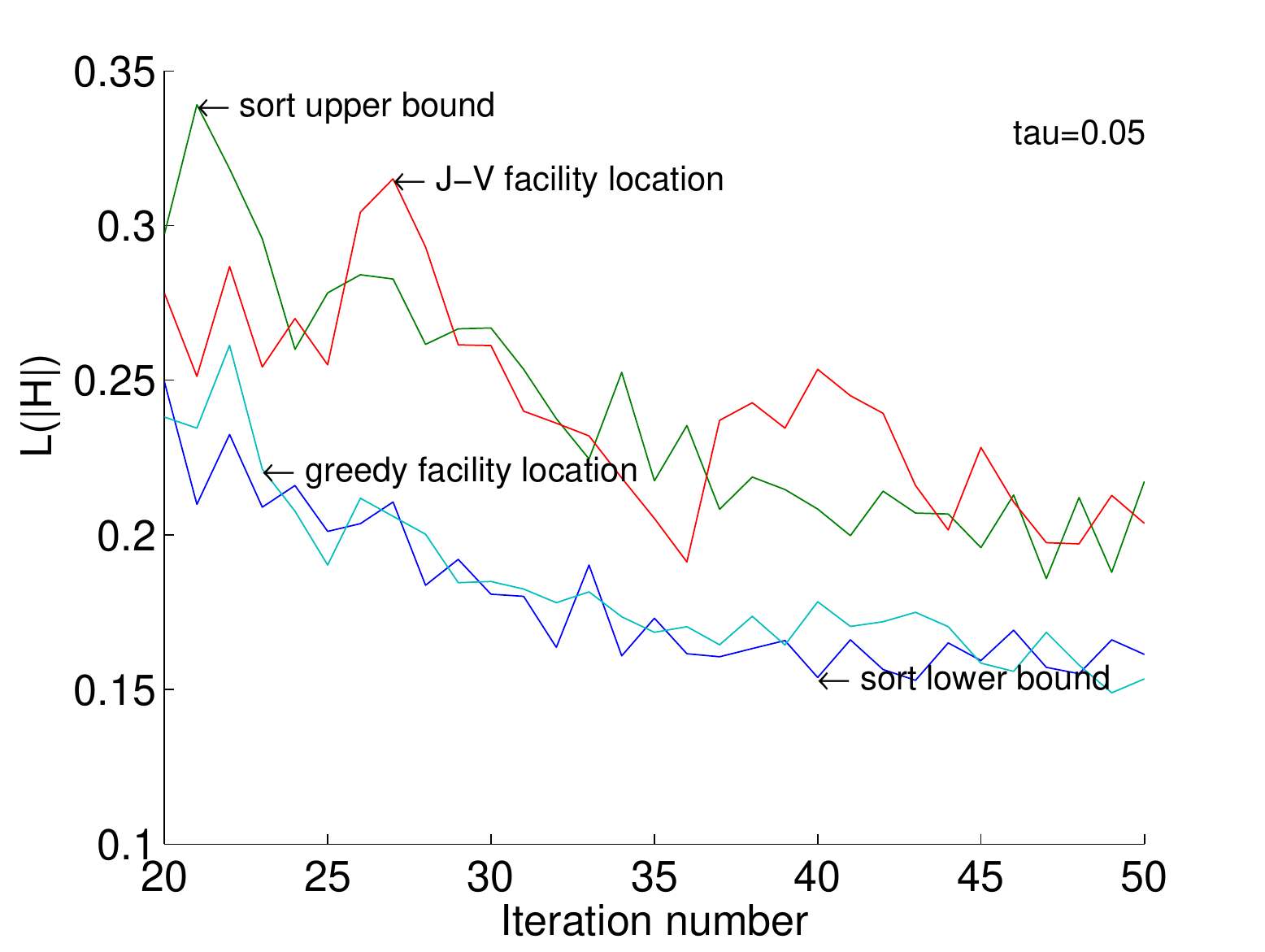}
      \caption{Discrete $L_1$ norm of $H$, 
$\deltat=0.1$ (top),  $\deltat=0.05$ (bottom)}
      \label{fl3}
   \end{figure}

Now we take smaller time-discretization step-sizes.
Figure~\ref{fl3} 
compares the four pruning methods with $\deltat=0.1$ and $\deltat=0.05$. They both show that the \textit{sort lower bound} and the \textit{greedy facility location} pruning method are better than the two others. 

\subsection{Discussion}
Our experimental results confirm that the total approximation error comes both from the approximation error of the Lax-Oleinik semi-group and from the pruning error. The error of approximation of the semi-group can be improved by decreasing the discretization-time step-size $\deltat$, while the pruning error depends on the pruning method. When $\deltat$ becomes small, the pruning appears to be the bottleneck. We introduced here new pruning methods, combining facility location algorithms and semidefinite relaxations,
which improve the final precision
(see Figure~\ref{fl3}
). However, these pruning methods remain
time-consuming (see Table~\ref{tiquette}), new ideas
are needed to develop more efficient methods.
Our experiments also show that the 
error is of order $O(\deltat)$, which is smaller than
the bound of $O(\sqrt{\tau})$ established in~\cite{mccomplex}.
This remains to be studied theoretically.

\begin{table}
\caption{\label{tiquette} CPU time}
\centering
\begin{tabular}{l|l|l|l|l}
   $\tau$=0.2, $K$=25 & Total time & Propagation & SDP & Pruning \\
\hline
   \emph{sort lower} &1.04h  & 1.85\%& 98.15\%& 0.00\%\\
\hline
   \emph{sort upper} &1.34h & 1.52\%&98.43\% & 0.05\% \\
\hline
   \emph{J-V p-d} & 1.38h & 1.45\%&89.47\%& 9.08\%\\
\hline
   \emph{greedy} & 1.43h & 1.63\%&97.84\%& 0.53\%\\
\end{tabular}
\end{table}
\begin{figure}[thpb]
      \centering
      \includegraphics[scale=0.4]{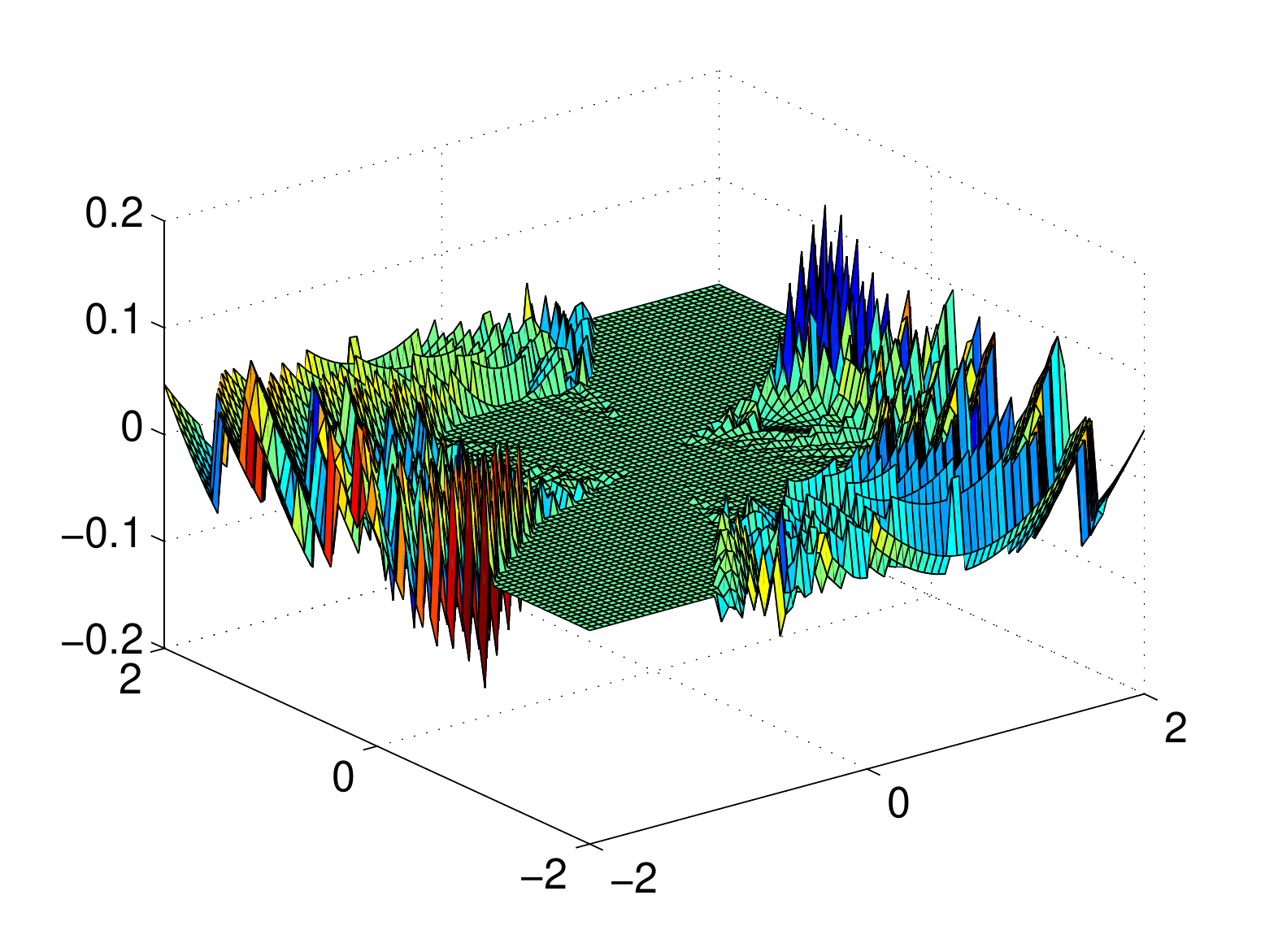}
      \includegraphics[scale=0.4]{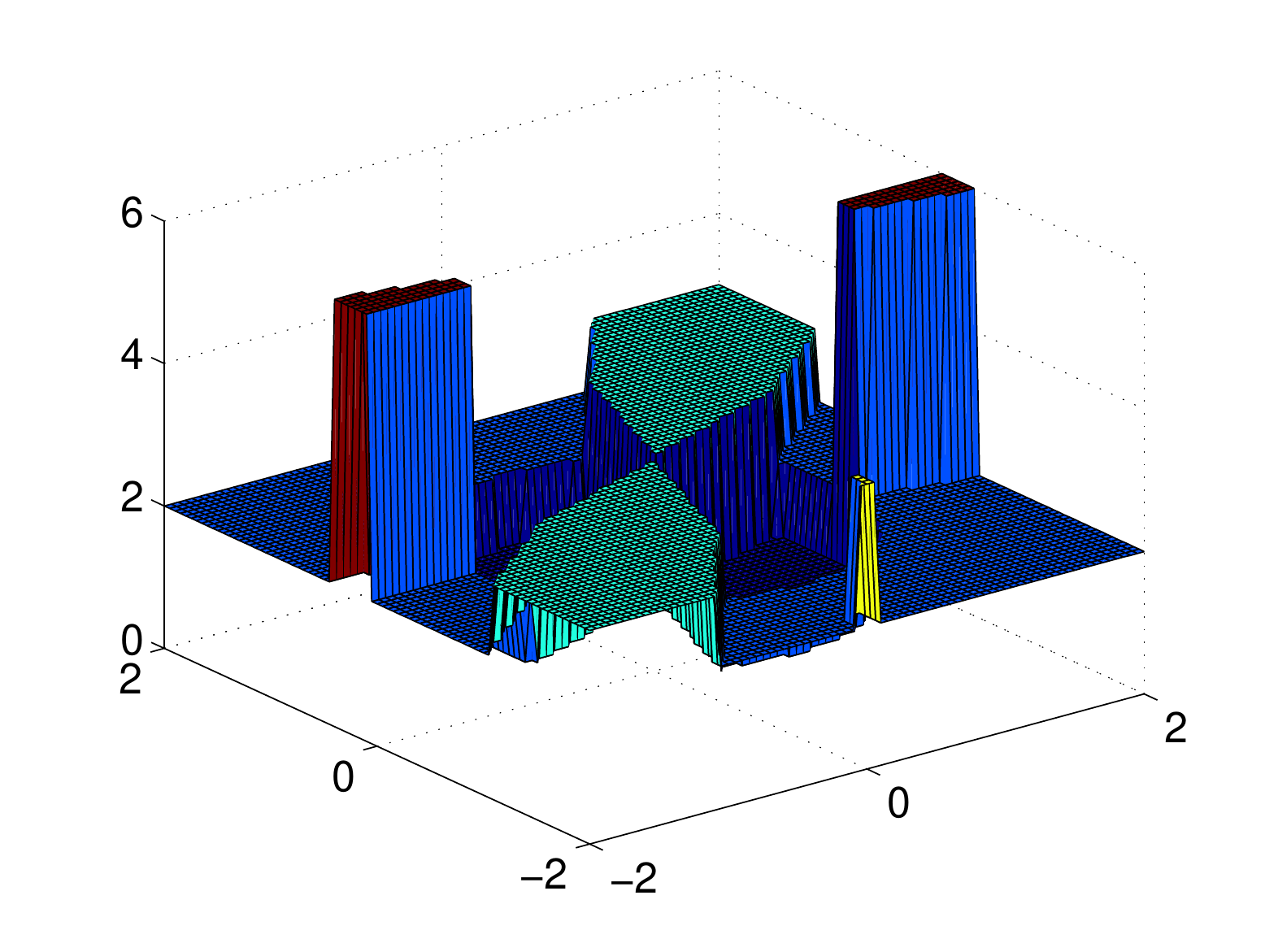}
      \caption{Backsubstitution error (top) and Optimal policy (bottom) on the $x_1$,$x_2$ plane, $\deltat=0.1$, with the \textit{greedy facility location} pruning method}
      \label{fl5}
   \end{figure}

\bibliographystyle{alpha}
\bibliography{biblio}

\begin{thebibliography}{MDG08}

\bibitem[AGL08]{a6}
M.~Akian, S.~Gaubert, and A.~Lakhoua.
\newblock The max-plus finite element method for solving deterministic optimal
  control problems: basic properties and convergence analysis.
\newblock {\em SIAM J. Control Optim.}, 47(2):817--848, 2008.

\bibitem[Br{\`e}67]{MR0215617}
L.~M. Br{\`e}gman.
\newblock A relaxation method of finding a common point of convex sets and its
  application to the solution of problems in convex programming.
\newblock {\em \u Z. Vy\v cisl. Mat. i Mat. Fiz.}, 7:620--631, 1967.

\bibitem[BV04]{MR2061575}
S.~Boyd and L.~Vandenberghe.
\newblock {\em Convex optimization}.
\newblock Cambridge University Press, Cambridge, 2004.

\bibitem[BZ07]{zidani-bokanowski}
O.\ Bokanowski and H.~Zidani.
\newblock Anti-dissipative schemes for advection and application to
  {H}amilton-{J}acobi-{B}ellman equations.
\newblock {\em J. Sci. Compt}, 30(1):1--33, 2007.

\bibitem[CD83]{capuzzodolcetta}
I.~Capuzzo~Dolcetta.
\newblock On a discrete approximation of the {H}amilton-{J}acobi equation of
  dynamic programming.
\newblock {\em Appl. Math. Optim.}, 10(4):367--377, 1983.

\bibitem[CFF04]{carlini-falcone-ferretti}
E.~Carlini, M.~Falcone, and R.~Ferretti.
\newblock An efficient algorithm for {H}amilton-{J}acobi equations in high
  dimension.
\newblock {\em Comput. Vis. Sci.}, 7(1):15--29, 2004.

\bibitem[CGQ04]{ilade}
G.~Cohen, S.~Gaubert, and J-P. Quadrat.
\newblock Duality and separation theorem in idempotent semimodules.
\newblock {\em Linear Algebra and Appl.}, 379:395--422, 2004.

\bibitem[CL84]{crandall-lions}
M.~G. Crandall and P.-L. Lions.
\newblock Two approximations of solutions of {H}amilton-{J}acobi equations.
\newblock {\em Math. Comp.}, 43(167):1--19, 1984.

\bibitem[DD97]{MRDrezner}
T.~Drezner and Z.~Drezner.
\newblock Replacing continuous demand with discrete demand in a competitive
  location model.
\newblock {\em NRL}, 44(1):81--95, 1997.

\bibitem[Fal87]{falcone}
M.~Falcone.
\newblock A numerical approach to the infinite horizon problem of deterministic
  control theory.
\newblock {\em Appl. Math. Optim.}, 15(1):1--13, 1987.
\newblock Corrigenda in {\it Appl. Math. Optim.}, 23:213--214, 1991.

\bibitem[Fer00]{MR1801497}
E.~Feron.
\newblock Nonconvex quadratic programming, semidefinite relaxations and
  randomization algorithms in information and decision systems.
\newblock In {\em System theory: modeling, analysis and control ({C}ambridge,
  {MA}, 1999)}, volume 518 of {\em Kluwer Internat. Ser. Engrg. Comput. Sci.},
  pages 255--274. Boston, MA, 2000.

\bibitem[FF94]{falcone-ferretti}
M.~Falcone and R.~Ferretti.
\newblock Discrete time high-order schemes for viscosity solutions of
  {H}amilton-{J}acobi-{B}ellman equations.
\newblock {\em Numer. Math.}, 67(3):315--344, 1994.

\bibitem[FM00]{a5}
W.~H. Fleming and W.~M. McEneaney.
\newblock A max-plus-based algorithm for a {H}amilton-{J}acobi-{B}ellman
  equation of nonlinear filtering.
\newblock {\em SIAM J. Control Optim.}, 38(3):683--710, 2000.

\bibitem[Gru93]{a2}
P.~M. Gruber.
\newblock Asymptotic estimates for best and stepwise approximation of convex
  bodies.\romannumeral 1.
\newblock {\em Forum Math.}, 5(5):281--297, 1993.

\bibitem[Gru07]{b2}
P.~M. Gruber.
\newblock {\em Convex and discrete geometry}.
\newblock Springer, Berlin, 2007.

\bibitem[Hla49]{MR0030547}
E.~Hlawka.
\newblock Ausf{\"u}llung und {\"u}berdeckung konvexer {K}{\"o}rper durch
  konvexe {K}{\"o}rper.
\newblock {\em Monatsh. Math.}, 53:81--131, 1949.

\bibitem[Lak07]{theseasma}
A.~Lakhoua.
\newblock {\em M{\'e}thode des {\'e}l{\'e}ments finis max-plus pour la
  r{\'e}solution num{\'e}rique de probl{\`e}mes de commande optimale
  d{\'e}terministe}.
\newblock Th\`ese de doctorat, Universit\'e Pierre et Marie Curie (Paris 6) et
  Universit\'e de Tunis El Manar, 2007.

\bibitem[LMS01]{litvinov}
G.~L. Litvinov, V.~P. Maslov, and G.~B. Shpiz.
\newblock Idempotent functional analysis: an algebraic approach.
\newblock {\em Math. Notes}, 69(5):696--729, 2001.

\bibitem[LV92]{Lin:1992:EMP:129712.129787}
J.-H. Lin and J.~S. Vitter.
\newblock e-approximations with minimum packing constraint violation.
\newblock In {\em Proceedings of the twenty-fourth annual ACM symposium on
  Theory of computing}, STOC '92, pages 771--782, New York, NY, USA, 1992. ACM.

\bibitem[McE06]{mceneaney-livre}
W.~M. McEneaney.
\newblock {\em Max-plus methods for nonlinear control and estimation}.
\newblock Systems \& Control: Foundations \& Applications. Birkh\"auser Boston
  Inc., Boston, MA, 2006.

\bibitem[McE07]{curseofdim}
W.~M. McEneaney.
\newblock A curse-of-dimensionality-free numerical method for solution of
  certain {HJB} {PDE}s.
\newblock {\em SIAM J. Control Optim.}, 46(4):1239--1276, 2007.

\bibitem[MDG08]{a7}
W.~M. McEneaney, A.~Deshpande, and S.~Gaubert.
\newblock Curse-of-complexity attenuation in the curse-of-dimensionality-free
  method for {HJB} {PDEs}.
\newblock In {\em Proc. of the 2008 American Control Conference}, pages
  4684--4690, Seattle, Washington, USA, June 2008.

\bibitem[MK10]{mccomplex}
W.~M. McEneaney and L.~J. Kluberg.
\newblock Convergence rate for a curse-of-dimensionality-free method for a
  class of {HJB} {PDE}s.
\newblock {\em SIAM J. Control Optim.}, 48(5):3052--3079, 2009/10.

\bibitem[NWF78]{MR0503866}
G.~L. Nemhauser, L.~A. Wolsey, and M.~L. Fisher.
\newblock An analysis of approximations for maximizing submodular set
  functions. {I}.
\newblock {\em Math. Programming}, 14(3):265--294, 1978.

\bibitem[Rei72]{MR0357936}
W.~T. Reid.
\newblock {\em Riccati differential equations}.
\newblock Academic Press, New York, 1972.
\newblock Mathematics in Science and Engineering, Vol. 86.

\bibitem[Rog64]{b1}
C.~A. Rogers.
\newblock {\em Packing and covering}.
\newblock Cambridge Tracts in Mathematics and Mathematical Physics, No. 54.
  Cambridge University Press, New York, 1964.

\bibitem[SSM10]{eneaneyphys}
M.R.~James S.~Sridharan, M.~Gu and W.M. McEneaney.
\newblock A reduced complexity numerical method for optimal gate synthesis.
\newblock {\em Phys. Review A}, 82(042319), 2010.

\bibitem[Vaz01]{MR1851303}
V.~V. Vazirani.
\newblock {\em Approximation algorithms}.
\newblock Springer-Verlag, Berlin, 2001.

\end{thebibliography}

\end{document}